\documentclass[11pt,a4paper]{amsart}
\usepackage{amssymb}
\usepackage{amscd}
\usepackage[all,2cell]{xy}
\SelectTips{cm}{}
\calclayout

\theoremstyle{plain}
\newtheorem{thm}{Theorem}[section]
\newtheorem{cor}[thm]{Corollary}
\newtheorem{lem}[thm]{Lemma}
\newtheorem{prop}[thm]{Proposition}

\theoremstyle{definition}
\newtheorem{defn}[thm]{Definition}
\newtheorem{exm}[thm]{Example}

\theoremstyle{remark}
\newtheorem{rem}[thm]{Remark}

\numberwithin{equation}{section}

\newcommand\lto{\longrightarrow}
\renewcommand\dim{\operatorname{dim}}
\newcommand\coh{\operatorname{Coh}\nolimits}
\newcommand\Nrep{\operatorname{NRep}\nolimits}
\newcommand\nrep{\operatorname{nrep}\nolimits}
\newcommand\mo{\operatorname{mod}\nolimits}
\newcommand\Mo{\operatorname{Mod}\nolimits}
\newcommand\Proj{\operatorname{Proj}\nolimits}
\newcommand\rep{\operatorname{rep}\nolimits}
\newcommand\Rep{\operatorname{Rep}\nolimits}
\newcommand\id{\operatorname{Id}\nolimits}
\newcommand\sid{\operatorname{id}\nolimits}
\newcommand\Cone{\operatorname{Cone}\nolimits}

\newcommand\Hom{\operatorname{Hom}\nolimits}
\newcommand\HH{\operatorname{HH}\nolimits}
\newcommand\End{\operatorname{End}\nolimits}
\newcommand\Ker{\operatorname{Ker}\nolimits}
\newcommand\ext{\operatorname{Ext}\nolimits}
\newcommand\im{\operatorname{Im}\nolimits}
\newcommand\soc{\operatorname{Soc}\nolimits}

\newcommand\rk{\operatorname{rk}\nolimits}
\newcommand\chara{\operatorname{char}\nolimits}
\newcommand\proj{\operatorname{proj}\nolimits}

\newcommand\ma{\mathcal{A}}
\newcommand\mh{\mathcal{H}}
\newcommand\mc{\mathcal{C}}

\newcommand\mt{\mathcal{T}}
\newcommand\ms{\mathcal{S}}
\newcommand\mcp{\mathcal{P}}
\newcommand\mcu{\mathcal{U}}

\newcommand\mbx{\mathbb X}
\newcommand\mbz{\mathbb Z}

\newcommand\bfC{\mathbf C}
\newcommand\bfD{\mathbf D}
\newcommand\bfK{\mathbf K}

\newcommand\mfp{\mathfrak p}

\newcommand\sg{\mathrm{sg}}

\begin{document}

\begin{abstract}
We discuss some basic properties of the graded center of a triangulated
category and compute examples arising in representation theory of
finite dimensional algebras.
\end{abstract}

\title{On the center of a triangulated category}
\author{Henning Krause}
\address{Henning Krause\\ Institut f\"ur Mathematik\\
Universit\"at Paderborn\\ 33095 Paderborn\\ Germany.}
\email{hkrause@math.upb.de}

\author{Yu Ye}
\address{Yu Ye\\ Department of Mathematics\\ University of Science and
Technology of China\\ Hefei 230026, Anhui\\ PR China.}
\email{yeyu@ustc.edu.cn}

\thanks{The second named author gratefully acknowledges support by the
Alexander von Humboldt Stiftung. He is supported in part by the
National Natural Science Foundation of China (Grant No. 10501041).}

\maketitle
\section{Introduction}

The graded center of a triangulated category $\mt$ with suspension
functor $\Sigma$ is a $\mathbb Z$-graded ring. The degree $n$
component consists of all natural transformations from the identity
functor $\id$ to $\Sigma^n$ which commute modulo the sign $(-1)^n$
with $\Sigma$.  The graded center is the universal graded commutative
ring that acts on $\mt$. For instance, the Hochschild cohomology
$\HH^*(A)$ of an algebra $A$ acts on the derived category $\bfD(A)$ via
a morphism $\HH^*(A)\to Z^*(\bfD(A))$ into the graded center.

It seems that the first systematic use of the graded center appears in
work of Buchweitz and Flenner on the Hochschild cohomolgy of singular
spaces \cite{BF}. For related work of Lowen and van den Bergh in the
setting of differential graded categories we refer to \cite{lv}.
Blocks of finite groups and their modular representation theory
provide the context for recent work of Linckelmann on the graded
center of stable and derived categories \cite{li}. Closely related is
the study of cohomological support varieties which depends on the
appropriate choice of a graded commutative ring acting on a
triangulated category; see \cite{bik}.

In this article, we prove some structural results and provide
complete descriptions of the graded center for some small examples.
The article is organized as follows.

In \S2, it is shown that for any abelian category $\ma$ with
enough projective objects, there is an isomorphism of graded
commutative rings $Z^*(\bfD^b(\ma))\cong Z^*(\bfD^{b}(\Proj(\ma)))$.
Here $\Proj(\ma)$ denotes the full subcategory of $\ma$ consisting of all
projective objects and the isomorphism is given by restriction.

In \S3-4, we deal with derived categories of hereditary
categories. Note that for a hereditary category, the derived
category and the bounded derived category have the same graded center.
In  \S3, the category $\mo(R)$ of finitely generated
modules over a Dedekind domain $R$ is considered. We calculate
$Z^*(\bfD(\mo(R)))$ explicitly. As we show it relates closely to the
residue fields of all the maximal ideals of $R$.

In \S4, we consider the module category of a tame hereditary algebra
and the category of coherent sheaves on a weighted projective line of
non-negative Euler characteristic.  We compute the graded centers of
their bounded derived categories.  Note that our methods do not apply
to wild cases. For a weighted projective line of wild type, we only
get a subalgebra of the graded center.

In \S\S5-6, we describe the graded centers of $\bfD^b(\mo(k[x]/x^2))$ and
$\underline{\mo}(k[x]/x^n)$ for $n\ge2$ respectively.

\section{Morphisms between graded centers}
\begin{defn}
Let $\mt$ be a triangulated category and $\Sigma$ the suspension
functor of $\mt$. We define a $\mbz$-graded abelian group
$Z^*(\mt)=Z^*(\mt, \Sigma)$ as follows. For any $n\in\mbz$, let
$Z^n(\mt)$ denote the collection of all natural transformations
$\eta\colon\id\to\Sigma^n$ which satisfy
$\eta\Sigma=(-1)^n\Sigma\eta$. The composition of natural
transformations gives $Z^*(\mt)$ the structure of a graded commutative
ring, and we call it the graded center of $\mt$.  Graded commutative
here means that $\eta\zeta=(-1)^{mn}\zeta\eta$ for all $\eta\in
Z^n(\mt)$ and $\zeta\in Z^m(\mt)$.
\end{defn}

\begin{rem}
(1) The definition of the graded center $Z^*(\mt)$ makes sense for any
graded category, that is, an additive category equipped with an
autoequivalence. In particular, the choice of the exact triangles of
$\mt$ is not relevant for $Z^*(\mt)$.

(2) The degree zero part $Z^0(\mt)$ is a subring of the usual center
$Z(\mt)$ of $\mt$, which by definition consists of all natural
transformations from the identity functor to itself. Note that
$Z^0(\mt)=Z(\mt)$ if $\Sigma=\id$.

(3) The graded center $Z^*(\mt)$ need not be a set in
general. However, it will be a set when the category $\mt$ is small.

(4) For any object $M$ in $\mt$ we define the graded ring
$\ext_{\mt}^{*}(M,M)$ by setting
$$\ext_{\mt}^{n}(M,M)=\Hom_{\mt}(M,\Sigma^{n}M)$$ for any integer
$n$. By definition there is a canonical graded ring morphism
\[Z^{*}(\mt)\longrightarrow\ext_{\mt}^{*}(M,M)\] mapping a natural transformation
$\eta\colon\id\to\Sigma^n$ to the morphism $\eta_{M}\colon M\to \Sigma^n M$.

Following Rouquier \cite{ro} we set $\langle M \rangle_1$ to be the full additive
subcategory of $\mt$ which contains $M$ and is closed under finite direct
sums, summands and the action of $\Sigma$, and for $i\ge2$ we define inductively
$\langle M \rangle_i$ as the full additive subcategory of $\mt$ consisting of all
objects isomorphic to direct summands of objects $Z$ for which there exists an exact
triangle $X\to Y\to Z\to\Sigma X$ with $X\in \langle M \rangle_1$ and
$Y\in \langle M \rangle_{n-1}$.

Now suppose that $M$ is an object of $\mt$ with $\mt =\langle M
\rangle_{d+1}$ for some positive integer $d$.  We set $\mathcal{N}$ to
be the kernel of the canonical morphism
$Z^{*}(\mt)\to\ext_{\mt}^{*}(M,M)$. It can be shown in this case that
$\mathcal{N}$ is a nilpotent ideal satisfying $\mathcal{N}^{2^{d}}=0$;
see \cite{li} for a proof. In particular, $Z^*(\mt)$ is modulo
nilpotent elements a set.
\end{rem}

Let $F\colon \ms\to \mt$ be an exact functor between triangulated
categories. An obvious question to ask is when the functor $F$ induces
morphisms between $Z^*(\ms)$ and $Z^*(\mt)$.  Recently, Linckelmann gave
an affirmative answer for this question in the case that there exists
a functor $G\colon\mt\to \ms$ which is simultaneously left and right
adjoint to $F$ and satisfies some further compatibility conditions
\cite{li}. The answer for general $F$ seems to be not known. The
following proposition shows that in some very specific situation, for
instance when $F$ is fully faithful, we do obtain some morphisms
between the graded centers.

\begin{prop}
Let $\mt$ be a triangulated category and $\ms$ a full triangulated
subcategory.
\begin{enumerate}
\item The inclusion functor $i\colon\ms\to\mt$ induces a morphism of
 graded rings $$i^*\colon Z^*(\mt)\lto Z^*(\ms),$$ where $i^*(\eta)_X =
 \eta_X$ for any  $\eta\in Z^*(\mt)$ and $X\in \ms$.
\item The canonical functor $\pi\colon\mt\to\mt/\ms$ induces a
morphism of graded rings $$\pi_*\colon Z^*(\mt) \lto Z^*(\mt/\ms),$$ where
$\pi_*(\eta)_X = \pi(\eta_X)$ for any  $\eta\in Z^*(\mt)$ and $X\in \mt/\ms$.
\end{enumerate}
\end{prop}
\begin{proof}
The proof is routine. To check that $\pi_*$ is well defined, one uses
the fact that for any commutative diagram in $\mt$
\[\xymatrix{
X \ar @{<--}[r]^{s} \ar[d]^{\alpha} &Z \ar[r]^{f}
\ar[d]^{\beta}&Y\ar[d]^{\gamma}\\ X'\ar @{<--}[r]^{s'}
&Z'\ar[r]^{f'}&Y' }
\]
with $\Cone(s)$ and $\Cone(s')$ in $\ms$, we have $\gamma\circ
(f/s)= (f'/s')\circ\alpha$ in the quotient category $\mt/\ms$, where we
use $\dashrightarrow $ to denote the morphisms whose cones are in
$\ms$.
\end{proof}

Up to now, little seems to be known about the properties of the above
morphisms. For example, the question when $i^*$ and $\pi_*$ are
surjective or injective is of special interest to us. Also, one might
study the induced morphism of graded rings $$(i^*, \pi_*)\colon
Z^*(\mt)\lto Z^*(\ms)\times Z^*(\mt/\ms).$$

\begin{exm}
Let $\ms\amalg\mt$ denote the direct product of two triangulated
categories $\ms$ and $\mt$. We view $\ms$ as a thick subcategory of
$\ms\amalg\mt$ and the corresponding quotient is equivalent to $\mt$.
Then we have $Z^*(\mt\amalg\ms)\cong Z^*(\mt)\times Z^*(\ms)$ via the morphism
$(i^*, \pi_*)$.
\end{exm}

For the rest of this section we focus on homotopy categories and
derived categories. Firstly we introduce some basic notations and
conventions. Let $\ma$ be any additive category. We denote by
$\bfC(\ma)$ the category of chain complexes in $\ma$. Recall that a
chain complex in $\ma$ is a sequence of morphisms in $\ma$
\[X = (\xymatrix{
  \cdots \ar[r]  &X_n \ar[r]^{d_{n}^X}
  & X_{n-1} \ar[r]& \cdots})
\]
with $d_{n}^Xd_{n+1}^X=0$ for all $n\in\mathbb{Z}$. A morphism of
complexes is a chain map $f\colon X \to Y $ consisting of a family of
morphisms $f_n\colon X_n\to Y_n$ in $\ma$ with $n\in \mathbb{Z}$ such
that $f_n\circ d_{n+1}^X= d_{n+1}^Y\circ f_{n+1}$ for all $n$, that
is, the diagram
\[\xymatrix{
  \cdots \ar[r]&X_{n+1} \ar[r]^{d_{n+1}^X}\ar[d]_{f_{n+1}} &X_n
  \ar[r]^{d_{n}^X}\ar[d]_{f_n} & X_{n-1} \ar[r]\ar[d]_{f_{n-1}}&
  \cdots\\ \cdots \ar[r] &Y_{n+1} \ar[r]^{d_{n+1}^Y} &Y_n
  \ar[r]^{d_{n}^Y} & Y_{n-1} \ar[r]& \cdots }
\]
commutes. We denote by $\bfC^+(\ma)$ the full subcategory of
$\bfC(\ma)$ which consists of all bounded below complexes, that is,
the complexes $X$ with $X_n=0$ for $n\ll0$. Similarly, we denote by
$\bfC^{-}(\ma)$ and $\bfC^{b}(\ma)$ the full subcategory of bounded
above complexes and complexes bounded in both directions, respectively.

If moreover $\ma$ is abelian, then for any integer $n$ the $n$-th
homology group $H_n(X)$ is by defintion $\Ker(d_n^X)/\im(d_{n+1}^X)$,
and any morphism $f$ of complexes induces morphisms of homology groups
$H_n (f)\colon H_n(X)\to H_n(Y)$ for all $n\in\mathbb{Z}$.

The homotopy category $\bfK(\ma)$ has the same objects as $\bfC(\ma)$.
The morphisms are the equivalence classes of the morphisms in
$\bfC(\ma)$ modulo the null-homotopic morphisms, that is, those with
components of the form
\[d_{n+1}^Y\circ h_n + h_{n-1}\circ d_n^X\]
for some family of morphisms $h_n\colon X_n\to Y_{n+1}$ in $\ma$,
$n\in\mathbb{Z}$.

The suspension functor (or shift functor) $\Sigma$ of $\bfC(\ma)$ is
defined by $(\Sigma X )_n = X_{n-1}$, $d_n^{\Sigma X}=-d_{n-1}^{X}$ on
the objects and by $(\Sigma f)_n= f_{n-1}$ on any morphism $f$.
Clearly $\Sigma$ is not only an autoequivalence but also an
automorphism of $\bfC(\ma)$. Moreover, $\Sigma$ also induces an
automorphism of $\bfK(\ma)$ and $\bfK(\ma)$ admits a triangulated
structure with  suspension functor $\Sigma$.

Let $\bfD(\ma)$ denote the derived category of $\ma$, i.e., the
localization of $\bfK(\ma)$ with respect to the quasi-isomorphisms.
Note that $\bfD(\ma)$ is again a triangulated category with suspension
functor $\Sigma$. One defines $\bfK^{*}(\ma)$ and $\bfD^{*}(\ma)$ with
$* \in\{+, b, -\}$ in a similar way.

Now let $\ma$ be an abelian category with enough projective objects,
and $\mathcal P$ be the full subcategory consisting of all projective
objects. We denote by $\bfK^{+,b}(\mathcal P)$ the thick subcategory
of $\bfK^+(\mathcal P)$ which consists of bounded below complexes $X$
with $H_n(X)=0$ for almost all $n$. Clearly, we have
$\bfK^{b}(\mathcal P)\subseteq \bfK^{+,b}(\mathcal P)$. In some cases,
objects in $\bfK^{b}(\mathcal P)$ are also called perfect complexes.

It is known that the following composition of functors
$$ \xymatrix{ \bfK^{+,b}(\mathcal P) \ar @{^{(}->}[d]_{ } \ar
@{-->}[rr]^{\cong} & & \bfD^b(\ma) \ar @{^{(}->}[d]^{ } \\
\bfK^+(\mcp)\ar[r]^{ }\ar @/^1pc/[rr]^{\cong} &\bfK^+(\ma) \ar[r]^{ }
& \bfD^+(\ma) }
$$ induces equivalences $\bfK^+(\mcp)\stackrel{\sim}\to \bfD^+(\ma)$
and $\bfK^{+,b}(\mathcal P)\stackrel{\sim}\to \bfD^b(\ma)$ of
triangulated categories. The quotient category
$\bfD_\sg(\ma)=\bfD^b(\ma)/\bfK^{b}(\mathcal P)$ is called the
triangulated category of singularities of $\ma$, because it is an
invariant of the singularities provided that $\ma$ is the category of
sheaves on some variety.  We know that $\bfD_\sg(\ma)=0$ if and only
if all objects of $\ma$ have finite homological dimension. When
$\ma=\Mo(A)$ for some self-injective ring $A$, then $\bfD_\sg(\ma)$ is
equavalent to the stable module category $\underline{\Mo}(A)$ of $A$.

We are now in a position to state our main result of this section.

\begin{thm}\label{iso-perfect}
 Let $\ma$ be an abelian category with enough projective objects and
$\mcp$ the full subcategory consisting of all projective objects.
Then the embedding $\bfK^{b}(\mcp)\to \bfK^{+,b}(\mcp)$ induces an
isomorphism $Z^*(\bfK^b(\mathcal P))\stackrel{\sim}\to Z^*(\bfD^b(\ma))$
of graded commutative rings.
\end{thm}

To prove the theorem, we need some preparations.

For each $n\in\mbz$ the $n$-th truncation functor
$\iota^n\colon\bfC(\ma)\to\bfC(\ma)$ is defined for a complex $X$ by
$(\iota^n X)_i=X_i$ for $i\le n$ and $0$ for $i>n$, and $d_i^{\iota^n
X}=d_i^X$ for $i\le n$ and $0$ for $i>n$. Clearly, $\iota^n$ sends
$\bfC(\ma)$ to $\bfC^{-}(\ma)$ and $\bfC^{+}(\ma)$ to
$\bfC^{b}(\ma)$. Note that we have a natural morphism $i^n\colon
\iota^n X\to X$, and sometimes we use $i^n_X$ to emphasize $X$. We
have $(i^n)_s=\sid$ for $s\le n$ and $0$ for $s> n$. The following
lemma is crucial in the proof of the main theorem.

\begin{lem}\label{nullhomotopy} Let $X\in \bfC(\mcp)$,
$Y\in \bfC(\ma)$ and $f\colon X\to Y$ be a chain map with $H_n(Y)=0$
for $n> 0$. Then $f$ is null-homotopic if and only if the
composition $f\circ i^n\colon \iota^n X\to X\to Y$ is null-homotopic
for some $n\ge 0$.
\end{lem}

\begin{proof} One direction is clear since the null-homotopic morphisms
form an ideal. Conversely, suppose that $f\circ i^n$ is null-homotopic
for some $n\ge 0$. To show that $f$ is also null-homotopic, it
suffices to find a family $\{h_n\colon X_n\to Y_{n+1}\mid
n\in\mathbb{Z}\}$, such that $f_n=d_{n+1}^Y h_n + h_{n-1} d_n^X$ holds
for all $n$.  By applying the shift functor, one can assume without
loss of generality that $f\circ i^0$ is null-homotopic. Thus there
exists a family of morphisms in $\ma$, say $\{h_n\colon X_n\to
Y_{n+1},n\le 0\}$, such that $f_n=d_{n+1}^Y h_n + h_{n-1} d_n^X$ for
all $n\le 0$.

Since $f$ is a chain map, we have $d_1^Y  f_1=f_0  d_1^X = d_1^Y h_0
d_1^X+ h_{-1}  d_0^X  d_1^X$, and hence $d_1^Y (f_1-h_0 d_1^X)=0$,
which implies that $\im(f_1-h_0
d_1^X)\subseteq\Ker(d_1^Y)=\im(d_2^Y)$, the last equality holds
because $H_1(Y)=0$. Now $X_1$ is projective implies that $f_1-h_0
d_1^X$ factors through $d_2^Y$, i.e., there exists $h_1\colon X_1\to
Y_2$ such that $f_1=d_{2}^Y  h_1 + h_{0}  d_1^X$, thus we get the
required $h_1$. Now repeat the argument and the lemma follows.
\end{proof}

\begin{prop}\label{extend}
Let $t\in \mathbb{Z}$ and $\eta\colon \id\to \Sigma^t$ be a natural
transformation for the category $\bfK^b(\mcp)$. Then $\eta$ extends
uniquely to a natural transformation
${\tilde{\eta}}\colon\id\to\Sigma^t$ for the category
$\bfK^{+,b}(\mcp)$.
\end{prop}

\begin{proof} First we will construct a morphism ${\tilde{\eta}}_X\colon X\to \Sigma^t
X$ for any $X\in \bfK^{+,b}(\mcp)$. The idea is to use certain
approximations.

Since $\eta$ is a natural transformation for
$\bfK^b(\mcp)$, we have for each $n$ a morphism
$\bar\zeta^n=\eta_{\iota^n X}\colon \iota^n X\to \Sigma^t\iota^n X$.
Now we fix a chain map $\zeta^0\colon \iota^0 X\to \Sigma^t\iota^0
X$ which is a representative of $\bar\zeta^0$. We can construct
inductively the representatives $\zeta^n$ of $\bar\zeta^n$ for all
$n\ge 0$, such that $\zeta^{n+1}_i=\zeta^n_i$ for all $n\ge 0$ and
$i\le n$.

In fact, suppose that $\zeta^n$ has been constructed, and let $\xi$
be any representative of $\bar\zeta^{n+1}X$. Consider the morphism
$j\colon \iota^n X\to\iota^{n+1}X$ which is given by $j_m=\sid_{X_m}$
for all $m\le n$ and 0 otherwise. Since $\eta$ is a natural
transformation, the diagram
 \[\CD
   \iota^nX @>j>> \iota^{n+1}X \\
   @V \zeta^n VV @V \xi VV  \\
   \Sigma^t\iota^nX @>\Sigma^t j>> \Sigma^t\iota^{n+1}X
 \endCD
 \]
commutes in the category $\bfK^b(\mcp)$, i.e., $\delta := \xi\circ j-
\Sigma^t j\circ \zeta^n$ is null-homotopic. Explicitly,
$\delta_i=\xi_i-\eta_i$ for $i\le n$ and $\delta_i=0$ for
$i\ge{n+1}$.

Now there exists a family of morphisms $\{h_i\colon (\iota^nX)_i\to
(\Sigma^t\iota^{n+1}X)_{i+1}\mid i\in\mathbb{Z}\}$ with $h_i=0$ for
$i>n$, such that $\delta_i=d_{i+1}^{\Sigma^t\iota^{n+1}X}\circ h_i +
h_{i-1}\circ d_i^{\iota^nX}$. The family $\{h_i\}$ can be viewed as a
family of morphisms $\{h_i\colon (\iota^{n+1}X)_i\to
(\Sigma^t\iota^{n+1}X)_{i+1}\mid i\in\mathbb{Z}\}$, thus it gives a
null-homotopic morphism $\delta'\colon \iota^{n+1}X\to
\Sigma^t\iota^{n+1}X$, which satisfies $\delta'_i=\delta_i$ for all
$i\le n$. We are done by setting $\zeta^{n+1}=\xi-\delta'$.

Now we define ${\tilde{\eta}}_X$ by $({\tilde{\eta}}_X)_n=\zeta^0_n$
for $n\le 0$ and $({\tilde{\eta}}_X)_n=\zeta^n_n$ for $n>0$. We
claim that ${\tilde{\eta}}$ is a natural transformation from $\id$
to $\Sigma^t$ for the category $\bfK^{+,b}(\mcp)$.

Note that by construction, ${\tilde{\eta}}_X$ satisfies the
following condition: for any $n\ge0$, there exists a representative
$\zeta^n_X$ for $\eta_{\iota^nX}$, which is given by
$(\zeta_{X}^n)_i=({\tilde{\eta}}_X)_i$ for $i\le n$, and
$(\zeta^n_X)_i= 0$ for $i\ge n+1$. In other words,
$\tilde\eta_X\circ i^n_X=\Sigma^t i^n_X \circ \eta_{\iota^nX}$ as
chain maps for all $n\ge0$, where $i^n_X$ denote the natural
morphism from $\iota^nX$ to $X$ as before.

Now let $X, Y\in \bfK^{+,b}(\mcp)$ and $f\colon X\to Y$ be any chain
map. Assume that $\tilde\eta_X: X\to\Sigma^t X$ and $\tilde\eta_Y:
Y\to\Sigma^t Y$ are arbitrary chain maps with the property
$\tilde\eta_X\circ i^n_X=\Sigma^t i^n_X \circ \eta_{\iota^nX}$ and
$\tilde\eta_Y\circ i^n_Y=\Sigma^t i^n_Y \circ \eta_{\iota^nY}$ for
$n\ge0$. We will show that ${\tilde{\eta}}_Y\circ
f=\Sigma^tf\circ{\tilde{\eta}}_X$ in $\bfK^{+,b}(\mcp)$. Note that in
the cube below, the other five faces are commutative by the
construction of $\tilde\eta_X$, $\tilde\eta_Y$, $\iota^n_X$ and
$\iota^n_Y$.
\[ \xymatrix@!=0.5pc{ & X \ar[rr]^{\bar\eta_X}
\ar'[d][dd]^{f} & & \Sigma^tX \ar[dd]^{\Sigma^tf} \\ \iota^nX
\ar[ur]^{i^n_X}\ar[rr]^{\eta_{\iota^nX}}\ar[dd] & & \Sigma^t\iota^nX
\ar[ur]\ar[dd] \\ & Y \ar'[r][rr]^{\bar\eta_Y} & & \Sigma^tY \\
\iota^nY \ar[rr]^{\eta_{\iota^nY}}\ar[ur]^{i^n_Y} & & \Sigma^t\iota^nY
\ar[ur] }
\]
By Lemma~\ref{nullhomotopy} it suffices to show that
${\tilde{\eta}}_Y\circ f\circ
i^n_X=\Sigma^tf\circ{\tilde{\eta}}_X\circ i^n_X$ for some
sufficiently large $n$. This is equivalent to showing that
${\tilde{\eta}}_Y\circ i^n_Y\circ \iota^n f=\Sigma^tf\circ\Sigma^t
i^n_X \circ \eta_{\iota^nX}$. The left hand side is
$\Sigma^ti^n_Y\circ\eta_{\iota^nY}\circ\iota^nf$, and since
$\eta_{\iota^nY}\circ \iota^n f= \Sigma^t\iota^nf\circ
\eta_{\iota^nX}$ and $\Sigma^tf\circ\Sigma^ti^n_X=
\Sigma^ti^n_Y\circ \Sigma^t\iota^n f$, the equality holds.

Thus by fixing such $\tilde\eta_X$ for each $X$, we can extend $\eta$
to the category $\bfK^{+,b}(\mcp)$. For the uniqueness, we need only
to take $f=\sid_X$ in the above argument. This completes the proof.
\end{proof}

\begin{cor}\label{positivegrad} Let $\ma$ and $\mcp$ be as before.
Then $Z^t(\bfK^b(\mcp))=0$ for all $t<0$, and therefore $Z^*(\bfK^b(\mcp))$
and $Z^*(\bfD^b(\ma))$ are positively graded.
\end{cor}

\begin{proof}
Suppose that $\eta$ is a natural transformation from
$\id_{\bfK^b(\mcp)}$ to $\Sigma^t_{\bfK^b(\mcp)}$ for some $t<0$. We
prove that $\eta_X=0$ by using induction on the length of the support
of $X$, where the support of $X$ means the interval $[i,j]$, such that
$i$ and $j$ are respectively the minimal and maximal integer $m$ with
$X_m\ne0$. Without loss of generality, we may assume that $i=0$ and we
use induction on $j$.

In the case $j=0$, clearly $\Hom_{\bfK^b(\mcp)}(X, \Sigma^tX)=0$ for
$t<0$.  Suppose $\eta_X=0$ for all $j\le m$ and suppose $X= (\cdots\to
0\to X_{m+1}\to\cdots\to X_0\to 0\to\cdots)$. By the same argument as
in the proof of Proposition~\ref{extend}, there is a representative
$\zeta$ of $\eta_X$, such that $\zeta_i=0$ for all $i\le m$, and now
the assumption $t<0$ implies that $(\Sigma^tX)_{m+1}=0$, which forces
that $\zeta_{m+1}=0$, thus $\zeta=0$ and hence $\eta_X=0$.
\end{proof}

With the above preparations, we can now prove the main theorem.

\begin{proof} [{Proof of Theorem~\ref{iso-perfect}.}]
Fix $\eta\in Z^t(\bfK^{b}(\mcp))$. By Proposition~\ref{extend}, $\eta$
extends uniquely to a natural transformation ${\tilde{\eta}}\colon \id
\to \Sigma^t $ for the category $\bfK^{+,b}(\mcp)$, and clearly
$i^*({{\tilde{\eta}})}=\eta$, where $i^*$ is induced by the embedding
$i\colon\bfK^{b}(\mcp)\to \bfK^{+,b}(\mcp)$. By the same argument as
in the last part in the proof of Proposition~\ref{extend}, one can
show that ${\tilde{\eta}}\Sigma=(-1)^n\Sigma{\tilde{\eta}}$, which
implies that ${\tilde{\eta}}\in Z^n(\bfK^{+,b}(\mcp))$. This proves
the surjectivity of $i^*$. The injectivity of $i^*$ follows from the
uniqueness of the extension.
\end{proof}

\begin{rem} Suppose there are enough injective objects in $\ma$
and denote by $\mathcal I$ the full subcategory of injectives.  Then
we have $\bfD^{b}(\ma)\cong \bfK^{-,b}(\mathcal I)$ and the dual
version of the theorem says that there is an isomorphism of graded
centers $ Z^*(\bfK^{b}(\mathcal I))\cong Z^*(\bfD^{b}(\ma))$.
\end{rem}

\section{Finitely generated modules over Dedekind domains}

The following two sections are devoted to studying the graded center of
the derived category of some hereditary categories. We look at some
basic examples and use explicit calculations. First we discuss the
derived category of the category of finitely generated modules
$\mo(R)$ for any Dedekind domain $R$. We start with some
preparation.

Let $R$ be an arbitrary unitary ring and denote by $Z(R)$ the center
of $R$. Let $z\in Z(R)$ and $M\in \mo(R)$. Then we have a morphism
$l_z\in \Hom_{R}(M, M)$, which is given by $l_z(m)=z\cdot m$. This is
indeed a morphism of modules since $z$ is in the center of
$R$. Moreover, $l_z$ induces a natural transformation from the
identity functor to itself for $\mo(R)$ as well as for $\bfD^b(\mo(R))$.

Now let $\mh$ be a hereditary abelian category, that is,
$\ext_\mh^i(M, N)=0$ for any $M, N\in\mh$ and $i\ge 2$. Consider the
derived category of $\mh$ and observe that any object $X\in D(\mh)$ is
isomorphic to $\bigoplus_{i\in \mbz}\Sigma^i(H_i(X))$. Here, $\Sigma$
is the shift functor and $H_i(X)$ is viewed as a stalk complex
concentrated in degree zero. For a simple proof of this, see \cite[\S1]{kr}. We
have the following easy lemmas.

\begin{lem}\label{degree0centerring}
Left multiplication induces an injective ring homomorphism
$Z(R)\to Z^0(\bfD^b(\mo(R)))$. Moreover, if $R$ is hereditary, then
this is an isomorphism.
\end{lem}

\begin{proof} 
For a proof, we just use the fact that left multiplication gives an
isomorphism from $Z(R)$ to the usual center of $\mo(R)$, i.e., the ring of
natural transformations from the identity functor to itself, and that
$\mo(R)$ is a full subcategory of $\bfD^b(\mo(R))$.

Moreover, if $R$ is hereditary, then all objects of $\bfD^b(\mo(R))$
are of the form $\bigoplus_{i\in \mbz}\Sigma^i M_i$ with
$M_i\in\mo(R)$ viewed as a stalk complex concentrated in degree
zero. Now the lemma follows easily.
\end{proof}

\begin{rem}
Note that the morphism in the lemma need not be an isomorphism; see
\cite{ku} or \S\ref{se:kx/x2} below.
\end{rem}

\begin{lem} Let $\mh$ be an arbitrary hereditary category. Then
$Z^*(\bfD^b(\mh))$ is concentrated in degree $0$ and $1$. Moreover,
the inclusions $\bfD^b(\mh)\subseteq \bfD^-(\mh)\subseteq \bfD(\mh)$ induce
isomorphisms of graded centers $Z^*(\bfD^b(\mh))\cong Z^*(\bfD^-(\mh))\cong
Z^*(\bfD(\mh))$.
\end{lem}

\begin{proof}
We have $\Hom_{\bfD(\mh)}(M, \Sigma^m M)=\ext^m_{\mh}(M, M)=0$ for all
$M\in \mh$ and $m\ge 2$, since $\mh$ is hereditary. Thus there is no
nontrivial natural transformations from $\id$ to $\Sigma^m$ for the
category $\bfD^b(\mh)$ for $m\ge 2$, and the first part of the lemma
follows. The last assertion follows from the fact that any element in
the graded center $Z^*(\bfD(\mh))$ is uniquely determined by the
restriction to the stalk complexes. The minor difference between both
derived categories is that any object in $\bfD(\mh)$ is an infinite
direct sum of stalk complexes while objects in $\bfD^b(\mh)$ can
always be written as finite direct sums.  Similarly, we have
$Z^*(\bfD^+(\mh))\cong Z^*(\bfD(\mh))$.
\end{proof}

Due to the lemma, to study the graded center of the derived categories
of hereditary abelian categories, we need only to consider the bounded
ones.

Now suppose that $\mh=\mh_1\vee\mh_2$, where $\mh_1$ and $\mh_2$ are
full additive subcategory of $\mh$, and we use $\vee$ to indicate that
any object of $\mh$ is a direct sum of an object of $\mh_1$ and an
object of $\mh_2$, and $\Hom_\mh(M_2,M_1)= \ext_\mh^1(M_1, M_2)=0$ for
all $M_1\in\mh_1$ and $M_2\in\mh_2$. We set $\Sigma^*\mh_1$ to be the
minimal additive subcategory of $\bfD^b(\mh)$ which contains $\mh_1$
and is closed under $\Sigma$, in other words the subcategory
consisting of all complexes with homologies contained in $\mh_1$. Note
that $\Sigma^*\mh_1$ is not a triangulated subcategory in general.
This will happen if $\mh_1$ is a thick subcategory of $\mh$, i.e.,
$\mh_1$ is closed under extensions, kernels and cokernels, and in this
case, $\mh_1$ is also a hereditary abelian category and
$\Sigma^*\mh_1\cong \bfD^b(\mh_1)$. For a proof of this, one uses
again the fact that any object in $\bfD^b(\mh)$ is a direct sum of
stalk complexes. Since $\Sigma$ is an autoequivalence of
$\Sigma^*\mh_1$, we can also define the graded center of
$\Sigma^*\mh_1$ with respect to $\Sigma$, and denote it by
$Z^*(\Sigma^*\mh_1)$.

\begin{prop}\label{degree1-wedge}
Let $\mh=\mh_1\vee\mh_2$ be a hereditary abelian category.  Then the
restriction map induces an isomorphism of abelian groups
\[Z^1(\bfD^b(\mh))\cong Z^1(\Sigma^*\mh_1)\times Z^1(\Sigma^*\mh_2).\]
\end{prop}

\begin{proof}
We produce an inverse map.  First observe that any object in
$\bfD^b(\mh)$ can be written uniqueley as $X_1\oplus X_2$ with
$X_1\in\Sigma^*\mh_1$ and $X_2\in\Sigma^*\mh_2$. Let $\eta_1\colon
\id_{\Sigma^*\mh_1}\to \Sigma_{\Sigma^*\mh_1}$ and $\eta_2\colon
\id_{\Sigma^*\mh_2}\to \Sigma_{\Sigma^*\mh_2}$ be natural
transformations. Then we define $\eta\colon \id_{\bfD^b(\mh)}\to
\Sigma_{\bfD^b(\mh)}$ by setting $\eta_{X_1\oplus X_2}$ to be the map
$(\eta_1)_{X_1}\oplus(\eta_2)_{X_2}$. We will show that $\eta$ is indeed
a natural transformation. To this end we need  to check that for
any morphism $f\colon X\to Y$ in $\bfD^b(\mh)$, we have $\Sigma
f\circ\eta_X= \eta_Y\circ f$.

Since any object of $\bfD^b(\mh)$ can be uniquely written as
$\bigoplus_{i\in\mbz}\Sigma^i(M_1^i \oplus M_2^i)$ with $M_1^i\in
\mh_1$ and $M_2^i\in\mh_2$, we need only to check the above
compatibility for the morphisms of the form $f\colon \Sigma^iM _1\to
\Sigma^j M_2$ and $g\colon \Sigma^i M_2\to\Sigma^j M_1$ with
$M_1\in\mh_1$ and $M_2\in\mh_2$. We claim that $\Sigma
f\circ\eta_{M_1}= \eta_{M_2}\circ f=0$ and $\Sigma g\circ\eta_{M_2}=
\eta_{M_1}\circ g=0$. In fact, since $\mh$ is hereditary, both sides
will vanish unless $j= i+1$ or $j=i$.  If $j=i+1$, the equalities hold
since $\ext_\mh^2(M,N)=0$ for all $M, N\in \mh$. Otherwise, if $j=i$,
we have $g=0$ and $\ext_\mh^1(M_1, M_2)=0$. Now the assertion follows
easily, and this completes the proof.
\end{proof}

Now we can begin the study of $Z^*(\bfD^b(\mo(R)))$ for a Dedekind
domain $R$. A Dedekind domain is an integral domain such that each
ideal can be written as a finite product of prime ideals, or
equivalently, a noetherian integrally closed domain with Krull
dimension at most one. This name was given to such rings in honor of
R. Dedekind, who was one of the first to study such rings in the
1870s. The rings of algebraic integers of number fields provide an
important class of Dedekind domains, which play a crucial role in
algebraic number theory.

The assumption on the Krull dimension implies that each nonzero
prime ideal of $R$ is maximal, and that the category $\mh=\mo(R)$ is
hereditary and any object $M$ of $\mo(R)$ is a direct sum of a
torsion-free module and a torsion module. Any finitely generated
torsion-free module is projective, and any finitely generated
torsion module is a finite direct sum of cyclic modules. Moreover,
by the Chinese Remainder Theorem, each cyclic module is a finite
direct sum of modules of the form $R/\mfp^l$ with $\mfp$ a maximal
ideal of $R$ and $l\ge 1$.

We have a decomposition of categories $\mh= \mh_+\vee\mh_0$, where
$\mh_+$ denotes the full subcategory of $\mh$ consisting of all
projective modules and $\mh_0$ consists of all torsion modules. Note
that $\mh_0$ is a thick subcategory and hence hereditary. Moreover,
there exists an Auslander-Reiten translation $\tau$ in $\mh_0$,
which is by definition an autoequivalence of $\mh_0$ such that there
exists a natural isomorphism $D\ext_\mh^1(X, Y)\cong \Hom_\mh(Y,
\tau X)$. This identity is usually called Serre duality and implies
the existence of Auslander-Reiten sequences in $\mh_0$.

Let $\max(R)$ denote the set of all maximal ideals of $R$. Then all
the indecomposables in $\mh_0$ are given by $\{R/ \mfp^l\mid
\mfp\in\max(R), l\ge 1 \}$. Denote by $\mh_{\mfp}$ the subcategory
of $\mh_0$ consisting of all $\mfp$-torsion modules, i.e., the
additive category generated by $\{R/\mfp^l\mid l\ge 1\}$ which is
also abelian, hereditary and a thick subcategory of $\mh_0$. Note
that $\mh_\mfp$ is $\tau$-invariant and $\mh_0=\coprod_{\mfp\in
\max(R)} \mh_\mfp$.

One can show that $\mh_\mfp$ is equivalent to the subcategory of
$\mo(R_\mfp)$ which consists of torsion $R_\mfp$-modules, where
$R_\mfp$ is the localization of $R$ at $\mfp$; and $\mh_\mfp$ is
obviously a full subcategory of $\mo(\hat{R}_\mfp)$, where
$\hat{R}_\mfp$ is the completion of $R$ with respect to $\mfp$,
i.e., the inverse limit $\varprojlim R/\mfp^l$. Note that we have an
isomorphism $R/\mfp \cong R_\mfp/\mfp R_\mfp$ of residue fields, and
we denote it by $k_\mfp$.

Fix an element $x\in \mfp\setminus \mfp^2$. Since $R_\mfp$ is a
discrete valuation ring, the multiplication with $x$ gives us an
isomorphism $R/\mfp^l\cong \mfp/\mfp^{l+1}$ of $R$-modules for any
$l\ge1$. Thus we have in $\mh_\mfp$ Auslander-Reiten sequences
\[\sigma^1_\mfp\colon 0\to R/\mfp\to R/\mfp^2\to R/\mfp\to 0\]
and\[\sigma^l_\mfp\colon 0\to R/\mfp^{l}\to R/\mfp^{l-1}\oplus
R/\mfp^{l+1} \to R/\mfp^{l}\to 0\] for all $l\ge 2$ , where the
morphism from $R/\mfp^{l+1}$ to $R/\mfp^l$ is the natural quotient
map, and we use the isomorphism $R/\mfp^l\cong \mfp/\mfp^{l+1}$
induced by the multiplication with $x$. This says that the
Auslander-Reiten quiver of $\mh_\mfp$ is a tube of $\tau$-period 1.
For the Auslander-Reiten sequences for Dedekind domains, see also
\cite{al}, Example3.1. Note that we have a natural equivalence
$\tau\cong\id_{\mh_\mfp}$.

It is easy to show that  $\Omega(R/\mfp^s)\cong R_\mfp$ for any $s\ge1$,
where $\Omega\colon \mo(R_\mfp)\to \mo(R_\mfp)$ is the syzygy
functor. We have a presentation of $R/\mfp^s$
$$\begin{CD}
  0 @> >> R_\mfp @>l_{x^s}>> R_\mfp @> >> R/\mfp^s @> >> 0
\end{CD}, $$
where $l_{x^s}$ denotes the multiplication by $x^s$. Note that
we have an isomorphism $R/\mfp^s\cong R_\mfp/\mfp^sR_\mfp$. Now the
above exact sequence induces an epimorphism $$\Hom_R(R_\mfp, R/\mfp^i)\to
\ext_R^1(R/\mfp^s, R/\mfp^i),$$ and when $i\le s$ this induces an
isomorphism. Now it is easy to show that the Auslander-Reiten sequence
$\sigma_\mfp^s$ corresponds to the composition $R_\mfp\to R/\mfp\cong
\mfp^{s-1}/\mfp^s\hookrightarrow R/\mfp^s$.  Moreover, all the
Auslander-Reiten sequences with starting term $R/\mfp^s$ are given by
the morphisms of the form $R_\mfp\to \soc(R/\mfp^s)\hookrightarrow
R/\mfp^s$, which equals $\lambda\sigma_\mfp^s$ for some
$0\ne\lambda\in k_\mfp$.

On the other hand, we identify $\Hom_{\bfD^b(\mh)}(M, \Sigma N) =
\ext^1_{\mh}(M, N)$ for any abelian category $\mh$ and any objects $M,
N\in \mh$. Now we can write down the graded center of
$\bfD^b(\mh_\mfp)$ explicitly by using the notion of a trivial
extension.

Let $A$ be an arbitrary ring and $M$ an $A$-$A$-bimodule. The trivial
extension ring of $A$ by the bimodule $M$, denoted by $T(A,M)$, is
defined to be the ring whose additive ring is $A\oplus M$ with
multiplication given by
$$(a,m)\cdot(a',m')=(aa',am'+ma')$$ for all $a,a'\in A$ and $m, m'\in
M$. Note that $T(A,M)$ can be identified with the subring of the upper
triangular ring
\[\begin{pmatrix}
A & M \\ 0 & A
\end{pmatrix}\]
which consists of all the matrices with equal diagonal entries.

The trivial extension ring is a positively graded ring which is
concentrated in degree $0$ and $1$ with $T(A,M)^0= A$ and
$T(A,M)^1=M$. Conversely, let $A= A^0\oplus A^1$ be an arbitrary
positively graded ring which is concentrated in degree $0$ and $1$.
Then $A\cong T(A^0, A^1)$ as graded rings, where the $A^0$-bimodule
structure on $A^1$ is induced by the multiplication of $A$. If $A$ is
commutative and $M$ an $A$-module, one can also define the trivial
extension $T(A, M)$, where $M$ is viewed as an $A$-$A$-bimodule. Note
that in this case $T(A, M)$ is always graded commutative.

\begin{prop}\label{homogeneoustube} Let $R$ be a Dedekind domain,
$\mfp$ a maximal ideal of $R$ and $k_\mfp$ the residue field. Then
as a graded ring,
 $$Z^*(\bfD^b(\mh_\mfp))\cong T(\hat{R}_\mfp, \prod_{l\in\mbz, l\ge1}k_\mfp),$$
where $k_\mfp$ is viewed as a simple $\hat{R}$-module.
\end{prop}

\begin{proof} Note that elements in $\hat{R}_\mfp$ are by definition
sequences $q=(q_i)_{i\in\mathbb{Z}, i\ge0}$ with $q_i\in R/\mfp^i$ and
satisfying $\pi_{i,j}(q_i)=q_j$ for all $i>j$, where $\pi_{i, j}\colon
R/\mfp^i\to R/\mfp^j$ is the natural quotient map. Now the collection
of morphisms $\{l_{q_i}\colon R/\mfp^i\to R/\mfp^i, i\in \mathbb{N}\}$
determines uniquely an element in $Z^0(\bfD^b(\mh_\mfp))$, where
$l_{q_i}$ is given by multiplication with $q_i$, and it is easy to
show this correspondence gives a bijection between $\hat{R}_\mfp$ and
$Z^0(\bfD^b(\mh_\mfp))$, which means that $Z^0(\bfD^b(\mh_\mfp))\cong
\hat{R}_\mfp$.

Now we consider the degree 1 component of the graded center. For any
$l\in \mbz, l\ge1$, we define $\eta^l_\mfp\in Z^1(\bfD^b(\mh_\mfp))$
by setting $(\eta^l_\mfp)_{R/\mfp^s}=0$ for all $s\ne l$,
$(\eta^l_\mfp)_{R/\mfp^l}=\sigma^l_\mfp$ and
$(\eta^l_\mfp)_{\Sigma^i(R/\mfp^s)}=(-1)^i\Sigma^i(\eta^l_\mfp)_{R/\mfp^s}$
for all $i, s$. To show $\eta^l_\mfp\in Z^1(\bfD^b(\mh_\mfp))$, it
suffices to show that $\eta^l_\mfp$ is a natural transformation. For
this, one needs to check that for all $i, j, m, n$ and $f:
\Sigma^iR/\mfp^m\to \Sigma^jR/\mfp^n$, the equality $\Sigma f\circ
(\eta^l_\mfp)_{\Sigma^i R/\mfp^m}= (\eta^l_\mfp)_{\Sigma^j R/\mfp^n}
\circ f$ holds. This is clear since both sides of the above equality
vanish, unless $f$ is an isomorphism, where we use the fact that any
$\sigma^l_\mfp$ is given by an almost split sequence.

The argument above shows that if for each $l$, we fix an
Auslander-Reiten sequence, say $\lambda_l \eta^l_\mfp$ for some
$\lambda_l\in k_\mfp$, with starting term $R/\mfp^l$, then we obtain
an element $\Sigma_{l}\lambda_l\eta^l_{\mfp}$ in
$Z^1(\bfD^b(\mh_\mfp))$. The infinite product makes sense since when
applied to any object in $\bfD^b(\mh_\mfp)$ it becomes a finite sum.

Conversely, let $\eta\in Z^1(\bfD^b(\mh_\mfp))$. We claim that
$\eta_{R/\mfp^l}$ corresponds to either an Auslander-Reiten sequence
or is zero for all $l$. This can be done by induction on $l$. Clearly it
holds for $l=1$. Assume that it holds for all $l\le n-1$, we will show
that it is also true for $l=n$.

From the exact sequence $0\to \mfp^{n-1}/\mfp^n\to R/\mfp^n
\overset{\pi}\to R/\mfp^{n-1}\to 0$ we obtain an exact sequence
$$\ext_R^1(R/\mfp^n, R/\mfp)\to \ext_R^1(R/\mfp^n, R/\mfp^n) \overset{f}\to
\ext_R^1(R/\mfp^n, R/\mfp^{n-1})\to 0.
$$ It is easy to show that $\Ker(f)$ is one dimensional and spanned by
the Auslander-Reiten sequences with starting term $R/\mfp^n$. Since
$\eta$ is a natural transformation, the map $\pi\colon R/\mfp^n\to
R/\mfp^{n-1}$ yields $\Sigma\pi\circ \eta_{R/\mfp^n}=
\eta_{R/\mfp^{n-1}}\circ\pi=0 $, where the last equality holds since
$\eta_{R/\mfp^{n-1}}$ is given by some Auslander-Reiten sequence or
zero. Thus $\eta_{R/\mfp^n}$ corresponds to an Auslander-Reiten
sequence by the above argument. The proposition now follows.
\end{proof}

Next we combine Lemma~\ref{degree0centerring},
Proposition~\ref{degree1-wedge} and that $\mo(R)= \mh_+\vee\mh_0$ with
$\mh_0= \coprod_{\mfp\in\max(R)}\mh_\mfp$. This gives the following
proposition. Note that $\mh_+$ consists of free modules and hence
$Z^1(\Sigma^*\mh_+)=0$.

\begin{prop} Let $R$ be a Dedekind domain and $\max(R)$ the set of
all maximal ideals. Then as a graded ring
$$Z^*(\bfD^b(\mo(R)))\cong T(R, \prod_{\mfp\in \max(R)}\prod_{l\in\mbz,
l\ge1}k_\mfp),$$ where each $k_\mfp\cong R/\mfp$ is viewed as a
simple $R$-module.
\end{prop}

\section{Tame hereditary algebras and weighted projective lines}

This section deals with the derived category for some further
classes of hereditary categories.  We consider either the category
of modules $\mo(A)$ of a tame hereditary algebra $A$ or the category
$\coh(\mbx)$ for a weighted projective line $\mbx$ of non-negative
Euler characteristic. Unfortunately, our methods do not work for the
wild cases. What we want to emphasize is that tubes are of special
importance in our calculations.

Throughout this section, $k$ denotes an algebraically closed field
and all categories considered are assumed to be $k$-linear; therefore
the graded centers are $k$-algebras. Note that most results hold for
an arbitrary base field $k$; however the proofs would require
modifications.

We begin by studying tubes. The tubes occurring in this section are
different from the ones for Dedekind domains and we will use
a different method to deal with them. Note that one can use
completed path algebras to unify the proofs.

Let $\mc$ be a uniserial hom-finite hereditary length
$k$-category. Recall that a length category is an abelian category
such that any object has a composition series of finite length. Note
that a length category is always a Krull-Remak-Schmidt category, i.e.,
any object can be written as a finite direct sum of indecomposables
and the endomorphism ring of any indecomposable object is
local. Following \cite{ar}, a locally finite abelian category is
called uniserial if any indecomposable object of finite length has a
unique composition series.

It follows from Theorem~2.13 in \cite{si}, that any hom-finite length
category is equivalent to the category of finite length comodules of
some basic coalgebra. And since $k$ is assumed to be algebraic closed,
any basic coalgebra is pointed, that is, it can be realized as a
subcoalgebra of certain path coalgebra of some quiver.

A quiver $Q =(Q_0, Q_1, s, t)$ is by definition an oriented graph,
where $Q_0$ is the set of vertices, $Q_1$ the set of edges which are
usually called arrows, $s$ and $t$ are two maps from $Q_1$ to $Q_0$
such that for each arrow $\alpha$, $s(\alpha)$ and $t(\alpha)$
denote respectively the starting vertex and the terminating vertex
of $\alpha$. A path in $Q$ is a sequence of arrows
$\alpha_1\alpha_2\cdots\alpha_n$ with $t(\alpha_i)=s(\alpha_{i+1})$
for $1\le i\le n-1$, $s(\alpha_1)$ and $t(\alpha_n)$ are called the
starting vertex and terminating vertex respectively and $n$ is the
length. Each vertex $v$ can be viewed as a path of length $0$ which
starts and terminates at $v$.

It is well known that there is a path algebra and a path coalgebra
structure on the vector space $kQ$ with basis consisting of all
paths in $Q$ and the multiplication and comultiplication are given by
composing and splitting the paths.  Denote by $kQ^a$ and
$(kQ^c,\Delta,\epsilon)$ the path algebra and path coalgebra of $Q$
respectively.

We denote the category of $k$-representations by $\Rep(Q)$ and the
subcategory of locally nilpotent representations by $\Nrep(Q)$. As
usual we denote the subcategories consisting of finite length objects
by $\rep(Q)$ and $\nrep(Q)$. It is well known that $\Rep(Q)$ is
equivalent to the module category of the path algebra $kQ^a$ and
$\Nrep(Q)$ is equivalent to the comodule category of the path
coalgebra $kQ^c$.

Let $n, m\in \mathbb{Z}\cup\{-\infty, +\infty\}$ with $n\le m$. We
use $A_{[n,m]}$ to denote the following quiver. The vertices are
indexed by $\{i\in \mathbb{Z}\mid n\le i\le m \}$ and for each $n\le
i\le m-1$ there is exactly one arrow which starts at the vertex $i$
and terminates at $i+1$. Now denote the quivers $A_{[-\infty, 0]},
A_{[-\infty, +\infty]}, A_{[0,+\infty]}$ and $A_{[0, n]}$ for any
$n\ge1$ by $A^{\infty}, A^{\infty}_{\infty}, A_{\infty}$ and $A_n$
respectively. Also we denote by $Z_n$ the basic cycle of length $n$
for any $n\ge1$, i.e., the quiver obtained from $A_n$ by gluing the
vertices $0$ and $n$.

The following classification is a special case of Theorem~2.10(i) in
\cite{cg}.

\begin{lem}\label{class-serial-cat} Let $\mc$ be a uniserial hereditary length $k$-category.
Then $\mc$ is equivalent to $\nrep(Q)$ for some quiver $Q$, where
$Q$ is a disjoint union of quivers of type $A^{\infty},
A^{\infty}_{\infty}, A_{\infty}, A_n$ or $Z_n$.
\end{lem}

The idea of the proof is easy. As shown in \cite{si}, $\mc$ is
equivalent to the category of finite length comodules of some pointed
coalgebra $C$, and that $\mc$ is hereditary means that $C$ must be a
path coalgebra and hence $\mc$ is given by some $\nrep(Q)$. The fact
that $\mc$ is uniserial implies that for each vertex $v\in Q_0$, there
is at most one arrow starting at $v$ and at most one arrow terminating
at $v$, and hence the lemma follows.

Now suppose that $Q$ is one of $A^{\infty}, A^{\infty}_{\infty},
A_{\infty}$ and $A_n$. Then the category $\nrep(Q)$ is directed. More
explicitly, any indecomposable object $M\in \nrep(Q)$ is a stone,
i.e., $\End_{\mc}(M)\cong k$ and $\ext^1_{\mc}(M,M)=0$. This has the
following consequence.

\begin{prop}\label{uniserial-directed}
$Z^*(\bfD^b(\nrep(Q)))\cong k$ for $Q= A^{\infty}, A^{\infty}_{\infty},
A_{\infty}$ and $A_n$.
\end{prop}

The only case left is $Q= Z_n$, where $n\ge1$ is a positive integer.
It is well known that $\nrep(Z_n)$ is a tube of $\tau$-period $n$.  We
need to fix some notation. Denote by $S_i$ the simple representation
with respect to the vertex $i$, and $M_i^{[l]}$ the indecomposable
representation with socle $S_i$ and of length $l$, for any $i$ and
$l\ge1$. In the category $\nrep(Q)$, there is neither a nonzero
projective nor a nonzero injective object, while in the category
$\Nrep(Q)$ there are enough injective objects, and we denote by
$M_i^{[\infty]}$ or simply $M_i$ the indecomposable injective module
with socle $S_i$. Note that $\{M_i^{[l]}\mid {0\le i\le n-1, 1\le
l\le\infty}\}$ gives a complete set of isoclasses of indecomposables
in $\Nrep(Q)$.

There is a monomorphism $i_s^l\colon M_s^{[l]}\to M_s^{[l+1]}$ and an
epimorphism $\pi_s^{[l]}\colon M_s^{[l]}\to M_{s-1}^{[l-1]}$ for any
$s$ and $l$, and any morphism between the indecomposables is a linear
combination of compositions of such morphisms.  For convenience, we
write again $\pi_s^{\infty}$ as $\pi_s$, and we set $M_s^{[l]}=0$ and
$i_s^l=\pi_s^l=0$ for $l\le 0$. More generally, we can define
monomorphisms
 \[i_s^{l,t}=i_s^{l+t-1}\circ\cdots\circ i_s^{l} \colon\ M_s^{[l]}\longrightarrow
 M_s^{[l+t]},\quad \forall\ 0\le s\le n-1, l\ge 1, t\ge 1,\]
and epimorphisms
 \[\pi_s^{l,t}=\pi_{s-t+1}^{l-t+1}\circ\cdots\circ\pi_s^l \colon\ M_s^{[l]}\longrightarrow
M_{s-t}^{[l-t]},\ \forall \ 0\le s\le n-1, l\ge 1, 1\le t\le l-1. \]
We also set $i_s^{l,\infty}$ to be the inclusion $M_s^{[l]}\to M_s$
and $i_s^{l,0}=\pi_s^{l,0}=\sid_{M_s^{[l]}}$. Note that we have the
equality of morphisms
 \[ \pi_{s}^{l+1}\circ i_s^{l}= i_{s-1}^{l-1}\circ\pi_s^{l}\colon\
M_s^{[l]}\longrightarrow M_{s-1}^{[l]}\] for all $s$ and $l$. The
syzygy functor $\Omega^{-1}$ is given by $\Omega^{-1}(M_s^{[l]})=
M_{s-l}$, $\Omega^{-1}(i_s^l)=\pi_{s-l}$ and
$\Omega^{-1}(\pi_s^l)=\sid_{M_{s-l}}$.  In the case $n=1$, the
subscript $s$ is omitted for simplicity.

\begin{lem}\label{degree0tube}
Let $Q= Z_n$. Then $Z^0(\bfD^b(\nrep(Q)))\cong k[[\xi]]$, where $\xi$ is
the natural transformation from the identity functor to itself, which
is given by
\[\xi_{M_s^{[l]}}=i_{s-n}^{l-n,n}\circ\pi_s^{l, n}\colon M_s^{[l]}\longrightarrow M_s^{[l]}. \]
\end{lem}

It is easy to check that $\xi$ is a natural transformation. Also the
infinite sum $\sum_{m\ge 0} \lambda_m \xi^m$ gives a natural
transformation, where $\lambda_m\in k$ for all $m$. Observe that this
does make sense, because the sum is indeed a finite sum when applied
to any object in $\nrep(Q)$. To show that this gives all the natural
transformations, one just uses the fact that $\{\xi_{M_s^{[l]}}^m, m\ge
0\}$ spans $\End_\mc(M_s^{[l]})$ for any $M_s^{[l]}\in \nrep(Q)$.

Clearly, we have an exact sequence $0\longrightarrow
M_s^{[l]}\overset{i_s^{l,\infty}}\longrightarrow
M_s\overset{\pi_s^{\infty,l} }\longrightarrow M_{s-l}\longrightarrow
0$ for any $M_s^{[l]}$. This induces an epimorphism
$\Hom_\mc(M_{r}^{[m]}, M_{s-l})\twoheadrightarrow
\ext_\mc^1(M_r^{[m]}, M_s^{[l]})$, which is an isomorphism when $m\le
l$. In particular, we can identify $\Hom_\mc(M_s^{[l]}, M_{s-l})$ with
$\ext_\mc^1(M_s^{[l]}, M_s^{[l]})$.

The following lemma is needed to describe the degree $1$ component of
the graded center of $\bfD^b(\nrep(Q))$.

\begin{lem}\label{dgree1-tubes}
Let $Q=Z_n$ and $n\ge 1$. If $n\ge 2$, then $Z^1(\bfD^b(\nrep(Q)))=0$; if
$n=1$, then as a $k$-vector space,
$$Z^1(\bfD^b(\nrep(Q)))\cong \prod_{l\in \mbz, l\ge1}k\cdot\eta^l,$$
where $\eta^l$ is given by $(\eta^l)_{M^{[l]}}= i^{1,\infty} \circ
\pi^{l, l-1}$ and $(\eta^l)_{M^{[a]}}= 0$ for $a\ne l$.
\end{lem}

\begin{proof}
First we consider the case $n\ge 2$. Fix $\eta\in
Z^1(\bfD^b(\nrep(Q)))$. We show that $\eta_{M_s^{[l]}}=0$ for all $s$
and $l$ by using induction on $l$. Clearly $\eta_{M_s^{[1]}}=0$ for
all $s$, since there is no self extension for the simple objects if
$n\ge 2$.

Now assume that the assertion holds for $l-1$. Applying the naturality of
$\eta$ to the injection $i_s^{l-1}\colon M_s^{[l-1]} \to M_s^{[l]}$,
we get the equality $\eta_{M_s^{[l]}}\circ i_s^{l-1} = \Sigma
i_s^{l-1}\circ \eta_{M_s^{[l-1]}}=0$. We claim that this equality holds only if
$\eta_{M_s^{[l]}}=0$. Otherwise, if $\eta_{M_s^{[l]}}\colon M_s^{[l]}\to
M_{s-l}$ is nonzero, then the dimension of $\im(\eta_{M_s^{[l]}})$ is at
least $n$ since $M_{s-l}^{[n]}$ is the minimal submodule of $M_{s-l}$
with the same top as $M_s^{[l]}$, thus the dimension of
$\im(\eta_{M_s^{[l]}}\circ i_s^{l-1})$ is at least $n-1$ and hence
nonzero (here we see the difference between $n=1$ case and $n\ge 2$ case), and now we use the
isomorphism $\Hom_\mc(M_{s}^{[l-1]}, M_{s-l})\cong
\ext_\mc^1(M_s^{[l-1]}, M_s^{[l]})$ to get that the left hand side of
the above equality is nonzero, this introduces a contradiction.

Now we assume that $n=1$. Note that any Auslander-Reiten sequence
with starting term $M_s^{[l]}$ is given by a nonzero multiple of
$\eta^l_{M^{[l]}}$. Now we can use the same argument as in the proof
of Proposition~\ref{homogeneoustube}.
\end{proof}

Combining Lemmas~\ref{degree0tube} and \ref{dgree1-tubes}, we get the
following.

\begin{prop}\label{center-basiccyle} Let $Q=Z_n$ and $n\ge1$.
If $n\ge 2$, then $Z^*(\bfD^b(\nrep(Q)))\cong k[[\xi]]$ is a graded
$k$-algebra concentrated in degree 0; if $n=1$, we have an
isomorphism
$$Z^*(\bfD^b(\nrep(Q)))\cong T(k[[\xi]], \prod_{l\in \mbz, l\ge1}k) $$
of graded algebras, where $k$ is viewed as the unique simple
$k[[\xi]]$-module on which $\xi$ acts trivially. Moreover, we have an
isomorphism of graded algebras $$T(k[[\xi]], \prod_{l\in \mbz,
l\ge1}k)\cong k[[\xi]][\eta]/(\eta^2),$$ where $\xi$ is of degree $0$ and
$\eta$ is of degree $1$.
\end{prop}

\begin{rem}
In case that $n=1$, we know that $\nrep(Z_1)$ is equivalent to
the category of finite dimensional nilpotent $k[x]$-modules, which is
just the category of finitely generated $(x)$-torsion modules over the
Dedekind domain $k[x]$. Thus Lemma~\ref{homogeneoustube} applies and
we get the same result. More generally, one can consider the completed
path algebra $k\hat{Z}_n$ of the quiver $Z_n$. Then $\nrep(Z_n)$ is
equivalent to the category of finite dimensional nilpotent modules
over $k\hat{Z}_n$, and the center of $k\hat{Z}_n$ is isomorphic to
$k[[x]]$.
\end{rem}

\begin{rem}
Combining Lemma~\ref{class-serial-cat} and
Propositions~\ref{uniserial-directed} and \ref{center-basiccyle}, we
have now a description of $Z^*(\bfD^b(\mc))$ for any uniserial
hereditary length $k$-categroy $\mc$.
\end{rem}

Next we consider the category of finite dimensional modules over
finite dimensional hereditary $k$-algebras. Since $k$ is assumed to be
algebraically closed, we need only consider the path algebras. Now let
$Q$ be a finite, connected quiver without oriented cycles and $A=kQ$
the path algebra. Note that in this case, the center of the algebra is
the field $k$. First we consider the finite type case.

\begin{prop}\label{finitetype}
Let $Q$ be a quiver such that the path algebra $kQ$ is of finite
representation type.  Then $Z^*(\bfD^b{\mo(kQ)})\cong k$.
\end{prop}
\begin{proof}
The proof is almost the same as the one of
Proposition~\ref{uniserial-directed}.  If $A$ is of finite
representation type, then any indecomposable $A$-module $M$ is a
stone. In particular, $\Hom_{\bfD^{b}(A)}(M, \Sigma^1
M)=\ext^{1}_{A}(M,M)=0$.
\end{proof}

Next we consider the tame case. Let $\tau$ be the Auslander-Reiten
translation in $\mo(A)$. The Auslander-Reiten quiver of $\mo(A)$
consists of the preprojective part, the preinjective part and the
regular part. Recall that a $A$-module $M$ is preprojective if and
only if $\tau^nM=0$ for sufficiently large $n$; and $M$ is
preinjecitve if and only if $\tau^{-n}M=0$ for sufficiently large
$n$. Modules without preprojective and preinjective summands are
called regular modules. Denote by $\mathcal P, \mathcal R$ and
$\mathcal I$ the full subcategory of preprojective modules, regular
modules and preinjective modules respectively. We have the
decomposition $\mo(A)=\mcp\vee\mathcal R\vee \mathcal I$. Let $\eta\in
Z^1(\bfD^b(\mo(A)))$, since preprojective and preinjective modules
have no self-extensions, we get the following easy lemma by applying
Proposition~\ref{degree1-wedge}.

\begin{lem}
Let $Q$ be a quiver such that the path algebra $kQ$ is of tame
representation type, and let $\mathcal R$ denote the full subcategory
of $\mo(kQ)$ consisting of regular modules. Then
$Z^1{\bfD^b(\mo(kQ))}\cong Z^1(\bfD^b(\mathcal R))$.
\end{lem}

Recall that for a tame quiver, the regular part of the
Auslander-Reiten quiver is a disjoint union of tubes, and there are
neither morphisms nor extensions between different tubes, i.e.,
$\mathcal R= \coprod_{\mathfrak{t}\in \mathfrak T} \mathcal
{R}_{\mathfrak{t}} $, where $\mathfrak T$ is an index set for all the
tubes, $ \mathcal{R}_{\mathfrak{t}}\cong \nrep(Z_{p(\mathfrak{t})})$
and $p(\mathfrak{t})$ denotes the $\tau$-period of $\mathcal R_\mathfrak{t}$.
Each tube is an abelian subcategory and we have $\bfD^b(\mathcal
R)\cong \coprod_{\mathfrak{t}\in \mathfrak{T}}
\bfD^b(\mathcal{R}_\mathfrak{t})$.  Applying
Propositions~\ref{center-basiccyle} and \ref{degree1-wedge}, we get the following.

\begin{prop}
Let $Q$ be a tame quiver, and $\mathfrak{T}_1$ the index set of
all homogeneous tubes. Then $$Z^*(\bfD^b(\mo(kQ)))\cong T( k,
\prod_{\mathfrak{t}\in\mathfrak{T}_1}\prod_{m\ge0}k).
$$
\end{prop}

Next we consider the weighted projective lines over the field $k$.
Recall that a weighted projective line $\mbx$ is defined through the
attached category $\coh(\mathbb X)$ of coherent sheaves, which is a
small $k$-category satisfying certain axioms. This concept was
introduced by Geigle and Lenzing in \cite{gl} to study the interaction
between preprojective modules and regular modules for tame
hereditary algebras. For a definition we refer the readers to
\cite[\S10]{le}, where one can also find most references about
this subject.

Weighted projective lines play an important role in the
classification of hereditary categories. By a theorem of Happel
\cite{ha}, any connected, Ext-finite, hereditary abelian
$k$-category which has a tilting complex is derived equivalent
either to the category $\mo(A)$ for some finite dimensional
hereditary algebra $A$ or the category $\coh(\mbx)$ for some
weighted projective line $\mbx$.

First we recall some basic facts. Let $\mbx$ be a weighted projective
line and $\mh=\coh(\mbx)$ the category of coherent sheaves. The
category $\mh$ has Serre duality, i.e., there exists an equivalence
$\tau\colon\mh\to\mh$ and a natural isomorphism $D\ext_\mh^1(X,
Y)\cong \Hom_\mh(Y, \tau X)$. We denote by $\mh_0$ the full subcategory
consisting of all objects of finite length. Then $\mh_0$ is a
hereditary abelian subcategory and $\mh_0=\coprod_{x\in C} \mcu_x$ for
some index set $C$, where $\mcu_x$ is a tube with finite $\tau$-period
$p(x)$. Members in $C$ are called the points of $\mh$. Note that there
are only finitely many points with $p(x)>1$.

 We denote by $\mh_+$ the subcategory consisting of all objects
without a simple subobject. Objects of $\mh_+$ are called vector
bundles. Any indecomposable object of $\mh$ is either of finite length
or a vector bundle. There is a linear form $\rk:
{\operatorname{K}_0}(\mh)\to\mbz$, called rank, which is
$\tau$-invariant, vanishes on objects of $\mh_0$ and takes positive
values on objects of $\mh_+$. Objects of $\mh_+$ of rank one are
called line bundles, and by definition $\mh$ contains a line
bundle. For any vector bundle $E$, we have a filtration $E_0\subseteq
E_1\subset\cdots\subseteq E_r =E$ with the line bundle factors
$E_i/E_{i-1}$, where $r=\rk(E)$.

For any line bundle $L$ and any point $x\in C$, $\sum_{S\in\mcu_x}
\dim_k\Hom_\mh(L, S)=1$, where $S$ runs through all simple objects in
$\mcu_x$. Clearly we have $\mh=\mh_+\vee\mh_0$, and therefore any
nonzero morphism between line bundles is a monomorphism.

Now we consider the graded center of $\bfD^b(\mh)$. Note that one can
define the Euler characteristic $\chi_\mh$ for $\mh$. If $\chi_\mh >
0$, then $\mh$ is derived equivalent to the category $\mo(A)$ for
some finite dimensional tame hereditary algebra $A$, and in this
case, the graded center has been computed. Firstly we have the
following easy lemma.
\begin{lem}
Let $\mbx$ be a weighted projective line. Then $Z^0(\bfD^b(\coh(\mbx)))=k$.
\end{lem}

\begin{proof}
We denote $\coh(\mbx)$ by $\mh$ as before. Since $\mh$ contains a line
bundle, we choose one and denote it by $L$.  Let $\eta\colon
\id_{\mh}\to \id_{\mh}$ be a natural transformation.  To prove the
lemma, it suffices to show that if $\eta_L=0$, then $\eta=0$.

Now assume that $\eta_L=0$. Let $x\in C$ be an arbitrary point, and
$S\in \mcu_x$ the simple object with $\Hom_\mh(L, S)\ne 0$. Note that
such $S$ exists and is unique. By Proposition~\ref{degree0tube},
$\eta_{\mcu_x}=0$ if and only if $\eta_{S^{[mr+1]}}=0$ for all $m\ge
0$, where $r=p(x)$ is the $\tau$-period and $S^{[mr+1]}$ is the
object in $\mcu_x$ with socle $S$ and of length $mr+1$. By using
induction we have $\dim_k\Hom_\mh(L, S^{[mr+1]})=m+1$ for any $m\ge0$.
We claim that there exists an epimorphism from $L$ to $S^{[mr+1]}$.
Otherwise, all morphisms will factor through $S^{[(m-1)r+1]}$, and
hence $\dim_k\Hom_\mh(L, S^{[mr+1]})=\dim_k\Hom_\mh(L, S^{[(m-1)r+1]})=m$,
which gives a contradiction.

Let $f\colon L\to S^{[mr+1]}$ be an epimorphism. Since $\eta $ is a
natural transformation, we have $\eta_{S^{[mr+1]}}\circ f=f\circ
\eta_L=0$, and hence $\eta_{S^{[mr+1]}}=0$. Now we have shown that
$\eta_{\mh_0}=0$. Conversely, using a similar argument, one can show
that if $\eta_{\mh_0}=0$, then $\eta_N=0$ for any line bundle
$N$. Since any vector bundle $E$ has a filtration with line bundle
factors, we get, using the five lemma, that $\eta_E=0$. This completes
the proof.
\end{proof}

Combined with Proposition~\ref{degree1-wedge},
Lemma~\ref{dgree1-tubes} and Proposition~\ref{center-basiccyle}, we
obtain the following embedding of algebras.

\begin{lem}
Let $\mbx$ be a weighted projective line, $\mh=\coh(\mbx)$ and $C_1$
the set of points of $\tau$-period $1$.  Then the algebra $Z=T(k,
\prod_{x\in C_1} \prod_{m\ge 0} k)$ is isomorphic to a subalgebra of
$Z^*(\bfD^b(\mh))$.
\end{lem}

In the tubular case, i.e., $\chi_\mh=0$, we have $\mh= \bigvee_{q\in
\mathbb Q\cup \{\infty\}} \mh^{\langle q\rangle}$, where for each $q$ we have
$\mh^{\langle q\rangle}\cong \mh_0$. In fact one can define the slope for
objects of $\mh$, and roughly speaking, for any $q\in \mathbb Q$,
$\mh^{\langle q\rangle}$ is just given by objects of slope $q$, and
$\mh^{\langle\infty\rangle}=\mh_0$. With this decomposition of categories, we
have the following proposition.

\begin{prop}
Let $\mbx$ be a weighted projective line of Euler characteristic 0, $\mh=\coh(\mbx)$ and
$C_1$ the set of points of $\tau$-period 1. Then
$$Z^*(\bfD^b(\mh))\cong T(k, \prod_{q\in \mathbb Q\cup \{\infty\}}
\prod_{x\in C_1}\prod_{m\ge 0} k). $$
\end{prop}

\section{The graded center of $\bfD^b(\mo(k[x]/(x^2)))$}
\label{se:kx/x2}

In this section, we will study the ring of dual numbers, which is by
definition the $k$-algebra $A=k[x]/(x^2)$, where $k$ is an arbitrary
base field. Set $\mc=\mo(A)$ and $\mathcal P$ the full subcategory of
$\mc$ consisting of projective modules. One has a complete description
of the indecomposable objects of $\bfD^b(\mc)=\bfK^{+,b}(\mathcal P)$,
and therefore one can write down the elements in $Z^*(\bfD^b(\mc))$
explicitly. By Theorem~\ref{iso-perfect}, we need only to consider the
category $\bfK^b(\mathcal P)$.

The indecomposable objects in $\bfK^{b}(\mathcal P)$ are well
understood, for example see \cite{ku}. They are given by
$\{A_m^n\mid-\infty < m\le n<\infty\}$, where $A_m^n$ is the complex
\[\entrymodifiers={!! <0pt, .8ex>+}
\xymatrix{\cdots \ar[r] &0 \ar[r] &{\underbrace{A}_n} \ar[r]^{x} & A
  \ar[r]^{x}& \cdots \ar[r]^{x}& {\underbrace{A}_m}\ar[r]
  &0\ar[r]&\cdots },
\]
that is, $(A_{m}^n)_i=A$ for $m\le i\le n$ and $0$ otherwise, and
$d_i^{A_{m}^n}= x$ for all $m< i\le n-1$, where we use $x$ to denote
the multiplication map $l_x$. If we allow $n$ to take the value
$\infty$, then we get all indecomposable objects in
$\bfK^{+,b}(\mathcal P)$, in fact $A^{\infty}_m\cong \Sigma^m S$ in
the derived category, where $S$ is the simple $A$-module, regarded
as a stalk complex concentrated in degree zero. Note that $\Sigma
A_m^n\cong A_{m+1}^{n+1}$. The following lemma is basic for our
computations.

\begin{lem}\label{basicmorphisms}
Let $-\infty < m\le n<\infty, -\infty < m'\le n'<\infty$. If $(m,n)\ne
(m',n')$, then $\Hom_{\bfK^{b}(\mathcal P)}(A_m^n, A_{m'}^{n'})$ is at
most one dimensional. The morphisms between indecomposable objects in
$\bfK^{b}(\mathcal P)$ are linear combinations of compositions of the
following four classes of morphisms:
\begin{enumerate}
\item[(a)] $\pi_m^{n,n'}\colon A_m^n\to A_{m}^{n'}$ for $m\le n'\le n$,
$(\pi_{m}^{n,n'})_m=x$ and $(\pi_{m}^{n,n'})_i=0$ $\forall i\ne m$;
\item[(b)] $\pi_{m,m'}^n\colon A_m^n\to A_{m'}^n$ for $m\le m'\le n$,
$(\pi_{m,m'}^n)_i=1$ $\forall m'\le i\le n$;
\item[(c)] $i_m^{n,n'}\colon A_m^n\to A_{m}^{n'}$ for $m\le n\le n'$,
$(i_{m}^{n,n'})_i=1$ $\forall m\le i\le n$;
\item[(d)] $i_{m,m'}^n\colon A_m^n\to A_{m'}^n$ for $m'\le m\le n$,
$(i_{m,m'}^n)_m=x$, and $(i_{m,m'}^n)_i=0$ $\forall i\ne m$.
\end{enumerate}
\end{lem}

The morphisms in the lemma look as follows.

\[(a)\qquad  \entrymodifiers={!! <0pt, .8ex>+}\xymatrix@!=0.4pc
 {0 \ar[r]  &{\overbrace{A}^n} \ar[r] &\cdots \ar[r] &{\overbrace{A}^{n'}} \ar[r]\ar[d]^{0} &\cdots
 \ar[r] & {\overbrace{A}^m}\ar[r]\ar[d]^{x} &0\\
 & & 0 \ar[r]  &A                 \ar[r] &\cdots \ar[r] &A                    \ar[r] &0}\]
\[(b)\qquad  \entrymodifiers={!! <0pt, .8ex>+}\xymatrix@!=0.4pc
 {0 \ar[r]  &{\overbrace{A}^n} \ar[r]\ar[d]^{1} &\cdots \ar[r] &{\overbrace{A}^{m'}} \ar[r]\ar[d]^{1} &\cdots
 \ar[r] & {\overbrace{A}^m}\ar[r] &0\\
  0 \ar[r]  &A                 \ar[r] &\cdots \ar[r] &A                    \ar[r] &0
 }
\]
\[(c)\qquad  \entrymodifiers={!! <0pt, .8ex>+}\xymatrix@!=0.4pc
 { & & 0 \ar[r]  &A\ar[d]^{1}                 \ar[r] &\cdots \ar[r] &A\ar[d]^{1}                    \ar[r] &0
 \\ 0 \ar[r]  &{\underbrace{A}_{n'}} \ar[r] &\cdots \ar[r] &{\underbrace{A}_{n}} \ar[r] &\cdots
 \ar[r] & {\underbrace{A}_m}\ar[r] &0}\]
\[ (d)\qquad  \entrymodifiers={!! <0pt, .8ex>+}\xymatrix@!=0.4pc
 {  0 \ar[r]  &A\ar[d]^{0}                 \ar[r] &\cdots \ar[r] &A\ar[d]^{x}                    \ar[r] &0
 \\ 0 \ar[r]  &{\underbrace{A}_n} \ar[r] &\cdots \ar[r] &{\underbrace{A}_{m}} \ar[r] &\cdots
 \ar[r] & {\underbrace{A}_{m'}}\ar[r] &0}
\]

The proof is straightforward and left to the reader. For any $m\le
n<\infty$, the space $\Hom_{\bfK^b(\mathcal P)}(A_m^n,A_m^n)$ is two
dimensional, and we denote the morphism $i_{m,m}^n=\pi_m^{n,n}$ by
$x_m^n$. Now let $\eta\colon \id_{\bfK^b(\mathcal P)}\to
\id_{\bfK^b(\mathcal P)}$ be a natural transformation. Clearly, $\eta$
is uniquely given by some datum $\{\mu_m^n, \lambda_m^n\in k,
-\infty<m\le n<\infty\}$ with $\eta_{A_m^n}=\mu_m^n\cdot 1+
\lambda_m^nx_m^n$.

\begin{prop}\label{degree0-dual-number}
Let $\eta\colon \id_{\bfK^b(\mathcal P)}\to \id_{\bfK^b(\mathcal P)}$
be a natural transformation and $\{\mu_m^n, \lambda_m^n\}$ the
corresponding datum.
\begin{enumerate}
\item We have $\mu_m^n=\mu_{m'}^{n'}$ for any $m, m', n$ and
$n'$. Conversely, any datum of the form $\{\mu, \lambda_m^n\in
k,\infty< m\le n<\infty\}$ arises as the datum of some natural
transformation $\eta\colon\id_{\bfK^b(\mathcal P)}\to
\id_{\bfK^b(\mathcal P)}$ by setting $\eta_{A_m^n}= \mu+
\lambda_m^nx_m^n$ for any $m$ and $n$.
\item If $\eta\in Z^{0}(\bfK^b(\mathcal P))$, then
$\lambda_m^n=\lambda_{m+r}^{n+r}$ for any $m, n$ and $r$, and any
elements in $Z^{0}(\bfK^b(\mathcal P))$ is obtained in this way.
\item As an algebra, $Z^{0}(\bfK^b(\mathcal P))\cong T(k, \prod_{r\ge0}
k)$, where $T(k, \prod_{r\ge0} k)$ is viewed as a graded
algebra concentrated in degree 0.
\end{enumerate}
\end{prop}
\begin{proof}
(1) Apply the naturality of $\eta$ to get $i_m^{n,n'}\circ
\eta_m^n=\eta_m^{n'}\circ i_m^{n,n'}$ and $\pi_{m,m'}^n\circ
\eta_m^n=\eta_m^{n'}\circ \pi_{m,m'}^{n}$. From this  follows that
$\mu_m^n=\mu_{m'}^{n'}$ for any $m, m', n$ and $n'$.

Conversely, for any datum of the form $\{\mu, \lambda_m^n\in k,\infty<
m\le n<\infty\}$, we claim that the above constructed $\eta$ is indeed
a natural transformation. In fact, one can easily show that the
equalities $f\circ\eta_m^n=\eta_{m'}^{n'}\circ f$ hold for those
morphisms $f\colon A_m^n\to A_{m'}^{n'}$ listed in
Lemma~\ref{basicmorphisms}.  Now using the fact that $\bfK^b({\mcp})$
is a Krull-Remak-Schmidt category and any morphism is some linear
combination of compositions of morphisms listed in
Lemma~\ref{basicmorphisms}, we get that $\eta_Y\circ f= f\circ \eta_X$
holds for any morphism $f\colon X\to Y$ in the category
$\bfK^b(\mcp)$. Thus $\eta$ is a natural transformation.

(2) Use the fact that by definition $\eta\in Z^{0}(\bfK^b(\mathcal
P))$ if and only if $\Sigma\eta=\eta\Sigma$, and this is equivalent to
the requirement that $\lambda_m^n=\lambda_{m+1}^{n+1}$ for any $m$ and
$n$.

(3) is an easy consequence of (2). In fact we can explicitly write
down the elements in $Z^{0}(\bfK^b(\mathcal P))$. For any $r\ge0$, let
$\eta_r\in Z^{0}(\bfK^b(\mathcal P))$ denote the natural transformation
obtained by setting $(\eta_r)_{A_m^n}=x_m^n$ for $m-n=r$ and 0
otherwise. Thus as vector spaces
$$Z^{0}(\bfK^b(\mathcal P))=k\cdot 1\oplus \prod_{r\ge0} k\cdot\eta_r. $$
By direct computation, the multiplication satisfies
$\eta_r\eta_{r'}=0$ for any $r$ and $r'$ and the isomorphism in (3)
follows.
\end{proof}

Now we consider the natural transformations from the identity functor
to $\Sigma^t$ for any positive integer $t>0$. Note that
$\Hom_{\bfK^b(\mathcal P)}(A_m^n,\Sigma^tA_m^n)=0$ for any $m, n$ with
$n<m+t$, and in the case $n\ge m+t$, the morphism space is one
dimensional with basis element $f_{t;m}^n = i_{m+t}^{n,
n+t}\circ\pi_{m, m+t}^{n}$. Let $\zeta\colon\id_{\bfD^b(\mathcal
P)}\to\Sigma^t$ be a natural transformation; it is uniquely determined
by the datum $\{\lambda_{t;m}^n, n\ge m+t\}$, where
$\zeta_{A_m^n}=\lambda_{t;m}^nf_{t;m}^n$. Applying the naturality of
$\zeta$ to the morphisms $i_m^{n,n'}$ and $\pi_{m,m'}^{n}$, one gets
$f_{t;m}^n=f_{t;m'}^{n'}$ for any $m, m', n$ and $n'$. Thus we get the
following lemma.

\begin{lem}
Let $t>0$. All natural transformations from $\id_{\bfK^b(\mathcal P)}$
to $\Sigma^t$ form a one dimensional $k$-space with a basis element
$\zeta_t$, where $\zeta_t$ is given by $(\eta_t)_{A_m^n}= f_{t;m}^n$
for all $n\ge m+t$ and $0$ otherwise. Moreover, the multiplication
satisfies $\zeta_t\zeta_{t'}=\zeta_{t+t'}$ for any $t,t'>0$ and
$\zeta_t\eta_r=\eta_r\zeta_t=0$ for any $t>0$ and $r\ge 0$, where the
$\eta_r$ are given as in the proof of
Proposition~\ref{degree0-dual-number}.
\end{lem}

Note that $\Sigma \zeta_t=(-1)^t\zeta_t\Sigma$ if and only if either
$\chara(k)=2$ or $\chara(k)\ne 2$ and $t$
is even. Combined with the last lemma and
Proposition~\ref{degree0-dual-number}, we get the graded center of
$\bfD^b(\mo(A))$.

\begin{prop}
Let $k$ be an arbitrary base field. Then as a graded algebra,
$$Z^*(\bfD^b(\mo(k[x]/(x^2))))\cong T(k[\zeta], \prod_{r\ge0}k),$$ where
$k$ is identified with $k[\zeta]/(\zeta)$ as a
$k[\zeta]$-module, $$Z^0(\bfD^b(\mo(k[x]/(x^2))))\cong T(k,
\prod_{r\ge0}k),$$ and $\zeta$ is of degree $2$ if
$\chara(k)\ne 2$, and of degree $1$ if
$\chara(k)=2$.
\end{prop}

\section{The graded center of the stable category $\underline{\mo}(k[x]/(x^n))$}

Another important class of triangulated categories are the stable
categories of self-injective algebras. We calculate the graded centers
in some special cases, namely for the algebras of the form
$k[x]/(x^n)$ with $n\ge 2$. These calculations are based on the fact
that the indecomposable objects and their morphisms are well
understood. Note that some algebras of the form $k[x]/(x^n)$ are
Brauer tree algebras, and we refer to \cite{keli} for the calculation of
the graded centers of their stable module categories.

Let $A=k[x]/(x^n)$ with $k$ an arbitrary base field. It is well known
that $A$ is uniserial and that all the indecomposable objects in
$\mo(A)$ are of the form $A_l=A/x^lA=x^{n-l}A$ with $1\le l\le
n$. There are epimorphisms $\pi_{r}^{l}=l_{x^{l-r}}\colon
A_l\twoheadrightarrow A_r$ for $l\ge r$ and monomorphisms
$i_{l}^r\colon A_l\hookrightarrow A_r$ for $l\le r$. For any $l$ and
$r$, $\Hom_A(A_l, A_r)$ has a basis $\{f^{l,r}_s=i_s^r\circ\pi_s^l\mid
1\le s\le \min(l,r)\}$. Moreover, the syzygy functor $\Omega$ is given
by $\Omega(A_l)=A_{n-l}$ and $\Omega(f^{l,r}_s)= \pi^{n-s}_{n-r}\circ
i^{n-s}_{n-l}=f^{n-l,n-r}_{n-r-l+s}$ for all $1\le l\le n-1$, $r,s \le
l$.

Now let $\mc=\underline{\mo}(A)$ be the stable category. One knows
that $\mc$ is a triangulated category with suspension functor $\Sigma
=\Omega^{-1}=\Omega$. In particular, we have $\Omega^2=\id_\mc$ in
$\mc$. The indecomposable objects in $\mc$ are given by
$A_l=A/x^lA=x^{n-l}A$ with $1\le l< n$, and $\bar f^{l,r}_s= 0$ if and
only if $l+r-n\ge s$. Consequently, $\Hom_\mc(A_l,
A_{n-l})=\Hom_A(A_l,A_{n-l})$.

For any self-injective ring $R$, let $Z(R)$ denote the graded center
of $R$.  There is a canonical morphism from $Z(R)$ to
$Z^0(\underline{\mo}(R))$. As we will show below, this map is not
injective in general. The more interesting question is that whether it
is surjective. In the case $A=k[x]/(x^n)$, the answer is yes. In fact,
for an arbitrary uniserial self-injective algebra, all natural
transformations from the identity functor to itself for the stable
category come from the center of the algebra.

\begin{prop}\label{degree0stablecat}
Let $A=k[x]/(x^n)$ with $n\ge2$ and $\mc=\underline{\mo}(A)$. Then
$Z^0(\mc)\cong k[x]/(x^{[\frac{n}{2}]})$, where $[\frac{n}{2}]$ denotes
the maximal integer which is no larger than $\frac{n}{2}$.
\end{prop}

\begin{proof}
Note that $A_{[\frac{n}{2}]}$ is of special
importance since $\End_\mc(A_{[\frac{n}{2}]})$ is of maximal
dimension among the indecomposable objects. Let $\eta$ be a natural
transformation from $\id_{\mc}$ to $\id_{\mc}$. We will show that
$\eta$ is uniquely determined by $\eta_{A_{[\frac{n}{2}]}}$.

We fix some $a\in A$ such that $\eta_{A_{[\frac{n}{2}]}}=\bar l_a$,
where $l_a$ is given by the multiplication with $a$ as before. Since
$\eta$ is a natural transformation, we have
$\pi_{[\frac{n}{2}]}^l\circ \eta_{A_l}=\bar
l_a\circ\pi_{[\frac{n}{2}]}^l$ for $l>[\frac{n}{2}]$ and
$i^{[\frac{n}{2}]}_l\circ \eta_{A_l}=\bar l_a\circ
i^{[\frac{n}{2}]}_l$ for any $l<[\frac{n}{2}]$. Now it is easy to show
that $\eta_{A_l}=\bar l_a$ for any $l$, since the solutions of the
equations above are unique. Therefore we have an epimorphism from $A$
to $Z^0(\mc)$, and easy computations show that $l_{x^{[\frac{n}{2}]}}=
0$ in $\mc$.
\end{proof}

Next we will compute the natural transformations from the identity
functor to $\Omega$. The following lemma is easy.

\begin{lem}\label{degree1stablecat}
Let $\zeta\colon \id_{\mc}\to \Omega$ be a natural
transformation. Then for any $1\le l<n$, we have $\zeta_{A_l}=
\lambda_l\cdot \bar f^{l,n-l}_1$ for some $\lambda_l\in k$. And
conversely, any family $\{\lambda_l, 1\le l< n\}$ induces a natural
transformation $\zeta$ by setting $\zeta_{A_l}= \lambda_l\cdot \bar
f^{l,n-l}_1$ for any $l$. Moreover, $\zeta\in Z^1(\mc)$ if and only if
$\lambda_l=-\lambda_{n-l}$ for any $l$.
\end{lem}

\begin{proof}
We use induction on $l$. Clearly we have $\zeta_{A_1}=\lambda_1\cdot
\bar f^{1,n-1}_1$. Now assume that $\zeta_{A_s}= \lambda_s\cdot \bar
f^{s,n-s}_1$ for some $\lambda_s\in k$, and consider the inclusion
$i_s^{s+1}$. One gets $\zeta_{A_{s+1}}\circ \bar
i_s^{s+1}=\bar\pi^{n-s}_{n-s-1}\circ \zeta_{A_s}=0$, and hence
$\zeta_{A_{s+1}}= \lambda_{s+1}\cdot \bar f^{s+1,n-s-1}_1$. The
remaining part is straightforward.
\end{proof}

Now let $\zeta_s$ denote the natural transformation given by
$(\zeta_s)_{A_l}= \delta_s^l\bar f^{s,n-s}_1$ for any $1\le l< n$.
We also denote by $t$ the identity map from $\id_\mc$ to
$\Omega^2=\id_\mc$ but viewed as an element in $Z^2(\mc)$.

Let $\tilde Z^*(\mc)$ be the $\mathbb Z$-graded space with $\tilde
Z^n(\mc)$ consisting of all natural transformations from $\id_\mc$ to
$\Omega^n$. Note that $\tilde Z^*(\mc)$ forms a graded algebra and
$Z^*(\mc)$ is a subalgebra of $\tilde Z^*(\mc)$.

Observe that the case $n=2$ is slightly different. In fact in this
case, not only  $\Omega^2$ but also the shift functor $\Omega$ itself
is equivalent to the identity functor. We deal with this case
separately. With the above notations, we get the following results.

\begin{prop}
Let $\mc=\underline {\mo}(k[x]/(x^2))$. Then $ \tilde
Z^*(\mc)=k[\zeta_1,\zeta_1^{-1}] $ with $\zeta_1$ of degree $1$. We have
$Z^*(\mc)= \tilde Z^*(\mc)$ if $\chara(k)=2$, and $Z^*(\mc) = k[\zeta_1^2]$
if $\chara(k)\neq 2$.
\end{prop}

Note that $\zeta_1^2$ equals $t$ as defined above, and clearly
$\zeta_1^{-1}$ is of degree $-1$. The proof follows directly
from Proposition~\ref{degree0stablecat} and
Lemma~\ref{degree1stablecat}.

\begin{prop}
Let $\mc=\underline {\mo}(k[x]/(x^n))$ and $n\ge 3$. Then we have
$$\tilde Z^*(\mc)= k[x, \zeta_1,\cdots,\zeta_{n-1}, t, t^{-1}]/\langle
x^{[\frac{n}{2}]}, x\zeta_s, \zeta_s x, \zeta_s\zeta_{s'}\rangle ,$$
where $x$, each  $\zeta_s$ and $t$ are in degree $0$, $1$ and $2$
respectively. Moreover, $Z^*(\mc)$ is the subalgebra generated by
$x$, $t$, $t^{-1}$, and $\zeta_s-\zeta_{n-s}$ with $1\le s\le
[\frac{n}{2}]$ if either $n$ is odd or $\chara(k)\ne 2$; if
$\chara(k)= 2$ and $n$ is even, then $Z^*(\mc)$ is the subalgebra
generated by $x$, $t$, $t^{-1}$, $\zeta_{[\frac{n}{2}]}$ and
$\zeta_s-\zeta_{n-s}$ with $1\le s\le [\frac{n}{2}]$.
\end{prop}

\begin{cor}
Let $\mc=\underline {\mo}(k[x]/(x^n))$ and $n\ge 3$. Then as a graded
algebra,
\[ Z^*(\mc)=k[x, \zeta_1,\cdots,\zeta_{l}, t, t^{-1}]/\langle
x^{[\frac{n}{2}]},
x\zeta_s, \zeta_s x, \zeta_s\zeta_{s'} \rangle
\] with $x$, each $\eta_s$ and $t$ in degree $0$, $1$ and $2$ respectively,
where $l= [\frac{n-1}{2}]$ if either $n$ is odd or $\chara(k)\ne 2$,
and $l=[\frac{n}{2}]$ if $\chara(k)=2$ and $n$ is even.
\end{cor}

\begin{rem}
For a self-injective algebra $A$, one has $\bfD^{b}(\mo(A))/
\bfK^{b}(\proj A)\cong \underline{\mo}(A)$. We have already seen
that $Z^*(\bfD^{b}(\mo(A)))\cong Z^*(\bfK^{b}(\proj A))$, but what can we say
about the ring homorphism $\pi_*\colon Z^*(\bfD^{b}(\mo(A)))\to
Z^*(\underline{\mo}(A))$?

For the algebra $A=k[x]/(x^2)$ we can describe $\pi_*$ explicitly, since
both graded centers are known. Recall that \[Z^*(\bfD^b(\mo(k[x]/(x^2))))=
(k\oplus \prod_{r\ge0}k\cdot\eta_r)[\zeta]/\langle \eta_r\eta_r',
\eta_r\zeta\rangle\] and $Z^*(\underline{\mo}(A))= k[t,t^{-1}]$. We know
that in this case, $\pi_*$ is neither injective nor surjective.
Explicitly, $\im(\pi_*)= k[t]$ and $\Ker(\pi_*)= \prod_{r\ge0}
k\cdot\eta_r$.
\end{rem}

\end{document}